\documentclass[11pt]{amsart}
\usepackage{mathrsfs, amsmath, amssymb, array,graphics, booktabs,diagbox,listings,color,textcomp,picinpar,wrapfig,cutwin,longtable,slashbox,amsthm,pst-node,tikz-cd}

\usepackage[page,title,titletoc,header]{appendix}
\usepackage[font=small,labelfont=bf]{caption}
\DeclareCaptionFont{tiny}{\tiny}
\DeclareCaptionFont{normalsize}{\normalsize}

\usepackage{geometry}
\definecolor{listinggray}{gray}{0.9}

\definecolor{lbcolor}{rgb}{0.9,0.9,0.9}

\usepackage{floatrow}
\newfloatcommand{capbtabbox}{table}[][\FBwidth]
\newcolumntype{P}[1]{>{\centering\arraybackslash}p{#1}}

\usepackage{algorithm}
\usepackage[noend]{algpseudocode}

\newtheorem{theorem}{Theorem}[section]
\newtheorem{lemma}[theorem]{Lemma}
\newtheorem{corollary}[theorem]{Corollary}
\newtheorem{proposition}[theorem]{Proposition}

\theoremstyle{definition}
\newtheorem{definition}[theorem]{Definition}

\theoremstyle{remark}

\numberwithin{equation}{section}

\begin{document}

\title{Algebraic fibrations of certain hyperbolic 4-manifolds}

\author{Jiming Ma}
\address{School of Mathematical Sciences \\Fudan University\\Shanghai 200433, China} \email{majiming@fudan.edu.cn}

\author{Fangting Zheng}
\address{Department of Mathematical Sciences \\ Xi'an Jiaotong Liverpool University\\ Suzhou 200433,	China  }
\email{Fangting.Zheng@xjtlu.edu.cn}

\keywords{algebraic fibered, hyperbolic 4-manifolds, finitely generated, finitely presented}
\subjclass[2010]{20F05, 20J05, 57M50, 57S25}
\date{Jan. 13, 2020}
\thanks{Jiming Ma was partially supported by  NSFC  11371094  and 11771088. Fangting Zheng was supported by XJTLU Research Development Funding RDF-19-01-29}

\begin{abstract}
	Algebraically fibering group is an algebraic generalization of the fibered 3-manifold group in higher dimensions. Let $M(\mathcal{P})$ and $M(\mathcal{E})$ be the cusped and compact hyperbolic real moment-angled manifolds associated to the hyperbolic right-angled 24-cell $\mathcal{P}$ and the hyperbolic right-angled 120-cell $\mathcal{E}$, respectively. Jankiewicz-Norin-Wise showed in \cite{JNW:2017}  that $\pi_1(M(\mathcal{P}))$ and $\pi_1(M(\mathcal{E}))$ are algebraic fibered. Namely, there are two exact sequences
	$$1\rightarrow H_{\mathcal{P}}\rightarrow \pi_1(M(\mathcal{P}))\xrightarrow{\phi_{\mathcal{P}}} \mathbb{Z}\rightarrow 1,$$
	$$1\rightarrow H_{\mathcal{E}}\rightarrow \pi_1(M(\mathcal{E}))\xrightarrow{\phi_{\mathcal{E}}} \mathbb{Z}\rightarrow 1,$$
	\noindent where $H_{\mathcal{P}}$ and $H_{\mathcal{E}}$ are finitely generated. In this paper, we furtherly show that the groups $H_{\mathcal{P}}$ and $H_{\mathcal{E}}$ are not $FP_2$. In particular, those fiber-kernel groups are finitely generated, but not finitely presented.
	\end{abstract}

\maketitle

\section{Introduction}

A group $G$ \emph{virtually algebraically fibers} if there is a finite index subgroup $G'$ admitting an epimorphism $G'\rightarrow \mathbb{Z}$ with finitely generated kernel $H$ \cite{JNW:2017}. The group $H$ is called the \emph{ fiber-kernel} of the algebraic fibration. Virtual fibrations of  3-manifold groups are fundamental problems in the study of 3-manifold topology. According to the well-known Stallings' criterion, provided $G$ is the fundamental group of an oriented, aspherical, compact 3-manifold $M$ and $G$ virtually algebraically fibers, then its fiber-kernel is the fundamental group of a surface $S$, and a finite-sheeted cover $M'$ of $M$, corresponding to a finite-index subgroup $G' < G$ , is an $S$-bundle over a circle \cite{Stallings:1962}.  Except a limited class of closed graph manifolds, every compact aspherical 3-manifold $M$ with $\chi(M)=0$ does virtually fiber \cite{Agol:2008,Agol:2013,Liu:2013,Wise:2018}.
 Kielak generalized the result of Agol and showed that a finitely generated virtually residually finite rational solvable (RFRS) group $G$ virtually admits algebraic fibration if and only if its first $L^2$-Betti number $\beta^{(2)}_1$ vanishes \cite{Kielak:2018:RFRS}. See also \cite{FV,Stover:2018} for virtually algebraic fibrations of  complex hyperbolic manifolds. Recently, Friedl and Vidussi use the BNS-invariant to show that the finitely generated group $G$ in an extension $1\rightarrow H\rightarrow G\rightarrow K\rightarrow 1$  algebraically fibers  if $b_1(G)>b_1(K)>0$ \cite{fv:2019}. Agol-Stover showed in \cite{as:2019} that certain Bianchi groups are virtually fibered on congruence subgroups .

Besides the finite generation, we want to say more about the finiteness properties of the fibered-kernel. In 1994, Bowditch and Mess \cite{Bowditch_Mess:1994} reported an example of a compact hyperbolic 4-manifold $M$, whose fundamental group $\pi_1(M)$ has a finitely generated subgroup $H$ while $H$ does not admit a finite presentation. Kapovich showed in \cite{Kapo:1998}  that, given a complex hyperbolic manifold $M$ of complex dimension of at least two, if $\pi_1(M)$ algebraically fibers, then the fiber-kernel is definitely not finitely presented.  See also \cite{Kapo:2013} for the non-coherence of arithmetic hyperbolic lattices and \cite{Krop:2018} for interesting result of fiber-kernels of hyperbolic groups. An aspherical manifold \emph{algebraically fibers} if its fundamental group algebraically fibers. In this paper, we consider two special four-dimensional real hyperbolic manifolds.
These are the first  twp examples of higher dimensional real hyperbolic manifolds that are  algebrically fibered with finitely generated, but not finitely presented fiber-kernels.
The finite generation was proved by \cite{JNW:2017}. Our contribution is to give a negative answer to the finite presentation of the fiber-kernel. The main theorem is as below:

\begin{theorem}\label{main}
	
Let $M(\mathcal{P})$ and $M(\mathcal{E})$ be the cusp and compact hyperbolic real moment-angled manifolds associated to the hyperbolic right-angled 24-cell $\mathcal{P}$ and the hyperbolic right-angled 120-cell $\mathcal{E}$, respectively. For the algebraic fibrations $$1\rightarrow H_{\mathcal{P}}\rightarrow \pi_1(M(\mathcal{P}))\xrightarrow{\phi_{\mathcal{P}}} \mathbb{Z}\rightarrow 1,$$
 $$1\rightarrow H_{\mathcal{E}}\rightarrow \pi_1(M(\mathcal{E}))\xrightarrow{\phi_{\mathcal{E}}} \mathbb{Z}\rightarrow 1$$
 \noindent given by Jankiewicz-Norin-Wise in \cite{JNW:2017}, the fiber-kernel groups $H_{\mathcal{P}}$ and $H_{\mathcal{E}}$ are not $FP_2$. In particular, they are finitely generated, but not finitely presented.
\end{theorem}

The outline of this paper is as follows. In Section 2, we give some preliminaries on the finiteness properties of groups, the right-angled Coxeter groups, the Bestvina-Brady's Morse theory on cubical complexes and the Jankiewicz-Norin-Wise's admissible systems. In Section 3 and Section 4, we prove Theorem \ref{main} in the cases of 24-cell and 120-cell, respectively.

\section{Virtually Algebraically Fibering of right-angled Coxeter group $C(\Gamma)$}
In this section, we review the Jankiewicz-Norin-Wise's legal system theory about virtually algebraic fibrations of right-angled Coxeter groups. The Bestvina-Brady's Morse theory matters a lot in this part. Some related terminologies about finiteness of a group are presented as well for the convenience of readers.

\subsection{Finiteness properties of groups}
We consider the following finiteness conditions for groups.
\begin{definition}
	A group $H$ is said to be of \emph{type $F_n$} if there is a Eilenberg-Mac Lane complex $K(H,1)$ with finite $n$-skeleton.
\end{definition}
Equivalently, a group is of type $F_n$ if it acts freely, properly, cellularly  and cocompactly on a ($n-1$)-connected cell complex. Clearly, a group is finitely generated if and only if it is of type $F_1$; and a group is finitely presented if and only if it is of type $F_2$.

By replacing the condition of ``($n-1$)-connected" with another condition of  ``integer-coefficient homological ($n-1$)-connected", we obtain a weaker finiteness property called $FH_n$. Clearly, if a group is of type $FH_1$, then it is of type $F_1$.

Another finiteness property $FP_n$, which is slightly weaker than the property $FH_n$, was introduced by Bieri \cite{Bieri:1976}:

\begin{definition}
	A group $H$ is said to be of \emph{type $FP_n$} if there exists a resolution$$P_n\rightarrow P_{n-1}\rightarrow\cdots\rightarrow P_0\rightarrow \mathbb{Z}$$of the trivial $\mathbb{Z}H$-module $\mathbb{Z}$ by finitely generated projective $\mathbb{Z}H$-modules $P_i$.
\end{definition}

By Schanuel's lemma, we have a standard way to show $H$ is of type $FP_n$ but not of type $FP_{n+1}$. Readers can refer to \cite{Brown:1982} for the detailed proof.

\begin{proposition} \label{Brown:1982} 
	If there is  a resolution
	$$Z_n\rightarrow P_n\rightarrow P_{n-1}\rightarrow\cdots\rightarrow P_0\rightarrow \mathbb{Z},$$
	where all the $P_i$ are finitely generated projective over $\mathbb{Z}H$, and $Z_n$ is not finitely generated over $\mathbb{Z}H$. Then $H$ is of type $FP_n$ but not of type $FP_{n+1}$.
\end{proposition}
Proposition \ref{Brown:1982} yields the following useful topological lemma. Then, we can use the topological information of the group action of $H$ on a cell complex to conclude the finiteness properties of $H$.
\begin{lemma} \label{B:Brown} \emph{(\cite{Brady:1999})}
	Suppose that the group $H$ acts freely, properly, cellularly, and cocompactly on the cell complex $X$. Then the following implications hold.
	
	(1) If $\widetilde{H_i}(X; \mathbb{Z})=0$ for $0\le i \le n-1$ and $\widetilde{H_n}(X; \mathbb{Z})$ is not finitely generated as a $\mathbb{Z}H$-module, then $H$ is of type $FP_n$ but not of type $FP_{n+1}$.
	
	(2) If X is ($n-1$)-connected, then $H$ is of type $FP_{n+1}$.
\end{lemma}

\subsection{\textbf{Right-angled Coxeter groups}}

\begin{definition}
	
 The \emph{right-angled Coxeter group} $C(\Gamma)$ associated to a finite simplicial graph $\Gamma$ is defined as	
 $$C(\Gamma)=\langle~ v_i\in V(\Gamma)~\vert~ v_i^2=1,~[v_i,v_j]=1~\text{if}~ (v_i,v_j)\in E(\Gamma)~\rangle,$$
where $V(\Gamma)$ and $E(\Gamma)$ are the vertex set and edge set of $\Gamma$, respectively.

\end{definition}

There is a canonical way to construct an Eilenberg-Mac Lane space for the group $C(\Gamma)$, for example see \cite{Davis:2008}. More precisely, the Coxeter group $C(\Gamma)$ acts properly and cocompactly on a CAT(0) cube complex $\widetilde{X}$, known as the \emph{Davis complex} of $C(\Gamma)$. The $1$-skeleton of $\widetilde X$ is isomorphic to the Cayley graph of $C(\Gamma)$ after identifying each bigon to an edge, and the $n$-cubes are equivariantly added to the $1$-skeleton for each collection of $n$ pairwise commuting generators. Note that there is a natural abelianization epimorphism $\alpha:C(\Gamma)\rightarrow \mathbb{Z}_2^{\vert V(\Gamma) \vert}$. Let $G'=\ker(\alpha)$ and $X = G'\backslash \widetilde X$. Jankiewicz-Norin-Wise \cite{JNW:2017} showed that such $X$ fits the setting of Bestvina-Brady's Morse theorey.

\subsection{Bestvina-Brady's Morse theory on cubical complexes}
Bestvina and Brady used the equivariant  Morse functions on cube complexes to study the finiteness properties of certain subgroups of right angled Coxeter groups \cite{BB:97,Brady:1999}. Here, we report some relevant definitions about the Bestvina-Brady's Morse theory on cubical complexes.

\begin{definition}\label{def:morse}
	Suppose that the group $G'$ acts freely, cocompactly, properly, by isometries on a contractible cubical complex $\widetilde{X}$ with $X=G'\setminus\widetilde{X}$ and thus, $\pi_1(X)=G'$. Let $\phi: G' \rightarrow \mathbb{Z}$ be an epimorphism, and let $\mathbb{Z}$ act on $\mathbb{R}$ by standard translations. The continuous map $\widetilde{\phi}:\widetilde{X}\rightarrow\mathbb{R}$ is called \emph{$\phi$-equivariant Morse function on $\widetilde{X}$} if:

	(1) $\widetilde{\phi}\circ g=\phi(g)\circ \widetilde{\phi} $ for all $g\in G'$.
	
	(2) For each $n$-cell $e$ of $\widetilde{X}$ with characteristic map $\chi_e:\square^n\rightarrow \widetilde{X}$, where $\chi_e(\square^n)=e$ and $\square^n$ is the standard $n$-cube, the composition $\widetilde{\phi}\circ \chi_e:\square^n\rightarrow \mathbb{R}$ extends to an affine map $\mathbb{R}^n\rightarrow \mathbb{R}$  and $\widetilde{\phi}\circ \chi_e$ is constant only for $n=0$.
	
	(3) The $\widetilde{\phi}$-image of the $0$-skeleton of $\widetilde{X}$ is discrete in $\mathbb{R}$.
	
\end{definition}

Note that if $\phi$ is induced by a map $\phi^{*}:X\rightarrow \mathbb{S}^1$, namely $\phi:\pi_1(X)\rightarrow \pi_1(\mathbb{S}^1)=\mathbb{Z}$ with kernel $H$,  it naturally satisfies the equivariant condition (1) in
Definition \ref{def:morse}. If the restriction of $\phi^{*}$ on every $n$-cube of $X$ is the composition of a
linear map $\mathbb{R}^n\rightarrow\mathbb{R}, (x_1,x_2,\cdots,x_n)\mapsto {\sum\limits_i} \pm x_i$ and the quotient
map $\mathbb{R}\rightarrow \mathbb{R}/\mathbb{Z}$, the map $\phi^{*}$ is called a \emph{diagonal map} \cite{JNW:2017}. The function $\widetilde{\phi}$ in diagram below would naturally satisfy the 

\[ \begin{tikzcd}
\widetilde{X} \arrow{r}{\widetilde{\phi}} \arrow[swap]{d}{} & \mathbb{R} \arrow{d}{} \\%
X \arrow{r}{\phi^*}& \mathbb{S}^1
\end{tikzcd}
\]

\noindent other two conditions (2) and (3) in Definition \ref{def:morse}.

It is observed in \cite{BB:97} that cubical complexes have a natural PL structure. Thus given a $n$-cell $e$ in a cubical complex $X$ with characteristic map $\chi_e:\square^n\rightarrow \mathbb{R}$, and a vertex $w\in \square^n$ mapping to $v\in e$, we have a well-defined map of links$$\chi_{e^*}:Lk(w,\square^n)\rightarrow Lk(v,e)\subset Lk(v,X).$$
The link $Lk(v,X)$ can be viewed as a union
$$Lk(v,X)=\bigcup\{\chi_{e^*}(Lk(w,\square^n))~\vert~ e~\text{is a cell of}~X~\text{and}~\chi_{e^*}(w)=v \}.$$

\begin{definition} 	Let $\widetilde{\phi}:\widetilde{X}\rightarrow \mathbb{R}$ be a $\phi$-equivariant Morse function. Define the \emph{ascending link} at a
vertex $w$ of the cube complex $X$ with respect to $\widetilde{\phi}$ to be
	$$A_w=\bigcup \{\chi_{e^*}(Lk(w,\square^n))~\vert~ \chi_{e^*}(w)=v ~\text{and}~ \widetilde{\phi}\circ\chi_e~\text{has a minimum at}~w~\},$$
	and the \emph{descending link} at $w$ with respect to $\widetilde{\phi}$ to be
	$$D_w=\bigcup \{\chi_{e^*}(link(w,\square^n))~\vert~ \chi_{e^*}(w)=v~ \text{and}~ \widetilde{\phi}\circ\chi_e~\text{has a maximum at}~w~\}.$$
\end{definition}

Bestvina and Brady showed in \cite{BB:97} and \cite{Brady:1999} the following theorem that simple topological conditions on links give rise to the finiteness properties for the kernel subgroup.

\begin{theorem}\label{B:finiteness}
	Let $\widetilde{\phi}:\widetilde{X}\rightarrow \mathbb{R}$ be a $\phi$-equivariant Morse function and let $H=Ker(\phi)$ be the kernel as defined above.
	
	(1) \emph{(Bestvina-Brady \cite{BB:97})}Suppose that each ascending link and each descending link is homologically $n$-connected. Then $H$ is of type $FH_{n+1}$.
	
	(2) \emph{(Brady \cite{Brady:1999})}Suppose the reduced homology of each ascending link and each descending link is zero in all dimensions $0$ through $n$+1, except for dimension $n$. Then $H$ is of type $FP_n$ but is not of type $FP_{n+1}$.
	
	
\end{theorem}
\subsection{\textbf{Jankiewicz-Norin-Wise admissible system}}
First, We introduce some useful combinatoric terminologies proposed by Jankiewicz-Norin-Wise \cite{JNW:2017} about  providing virtually algebraic fibrations. Let $\Gamma=\Gamma(V,E)$ be a simplicial graph with vertex set $V$ and edge set $E$.

A \emph{state} of $\Gamma$ is a subset $S\subset V$.
The state is \emph{legal} if the subgraphs induced by $S$ and by the complement $V-S$ of $S$ in $V$ are both nonempty and connected.
\begin{definition}\label{move}
A \emph{move at} $v\in V$ is an element $m_v\in 2^V$ with the following property:
\begin{enumerate}
	\item\label{move prop:1} $v\in m_v$.
	\item\label{move prop:2} $u\not\in m_v$ if $\{u,v\}\in E$.
\end{enumerate}
\end{definition}
A \emph{move system} is a choice $m_v$ of a move for each $v\in V$. We do not assume that $m_v\neq m_u$ for $v\neq u$.

Identify $\mathbb{Z}_2^{\vert V\vert}$ with $2^V$ in the obvious way where $\emptyset$ corresponds to the identity element
and multiplication in $\mathbb{Z}^{\vert V \vert}_2$ corresponds to the symmetric difference in $2^V$.
%
Both state and moves can be identified with some  elements of $\mathbb{Z}_2^{\vert V\vert}$.
A move system generates a subgroup $M$ of $\mathbb{Z}_2^{\vert V\vert}$.
The move system is \emph{legal} if there is a starting state $S$ such that all of elements in the $M$-orbit $M \cdot S$ are legal states.
We say such an orbit a \emph{legal orbit} and such a pair of legal move system and starting state an \emph{admissible system}.





By Theorem \ref{B:finiteness}, the connectness condition implies the finiteness properties for the kernel subgroup. Because the connectedness of a simplicial complex is determined by the connectedness of its 1-skeleton, we focus entirely on the 1-skeleton when discussing the ascending and descending links. Note that a diagonal map is determined by directing the 1-cubes of $X$ so that the orientations of opposite 1-cubes of each 2-cube agree. Jankiewicz-Norin-Wise showed the following combinatoric argument of legal system to provide virtually algebraic fibering of right-angled Coxeter group $C(\Gamma)$. The proof in \cite{JNW:2017} is quite sketchy, we report more details here for the convenience of readers.

\begin{theorem}\label{thm:win fiber} \emph{(Jankiewicz-Norin-Wise \cite{JNW:2017})}
	Let $\Gamma$ be a finite simplicial graph.
	Suppose there is an admissible system of legal move system and starting state over the graph $\Gamma$, then there is a
	Bestvina-Brady's Morse function $f:\widetilde{X}\rightarrow \mathbb{R}$ whose ascending and descending links of the $0$-cube of $\widetilde{X}$
	are non-empty and connected, where $\widetilde{X}$ is the cubical complex corresponding to $C(\Gamma)$, $X=G^{'}/\widetilde{X}$, $V$ is the vertex set of $\Gamma$ and $G^{'}=\text{ker}~(C(\Gamma)\rightarrow\mathbb{Z}_2^{\vert V \vert})$.
\end{theorem}
\begin{proof}
	Let $S$ be a state such that each element in $\langle m_v : v\in V\rangle S$ is legal. 
	Identify $z_i \in \mathbb{Z}_2^{\vert V\vert }$ with a set of vertices $V_i\in 2^V$ in the obvious way  where $\emptyset$ is the identity element 	and multiplication is the symmetric difference. Namely we have a map $\varphi: \mathbb{Z}_2^{\vert V\vert } \rightarrow 2^V,z_i\mapsto V_i\subset V $
	
	Consider a base $0$-cube $x\in X^0$. For each $ z_i \in \mathbb{Z}_2^{\vert V\vert }$, we  label other $0$-cubes of $X$ as $z_ix$. Suppose the $0$-cubes that connected to $z_ix$ by $1$-cubes are denoted by $z_{i_{1}}x,z_{i_{2}}x,\cdots,z_{i_{2^{\vert V \vert}}}x$, the set of the differences of  $z_i$
 with $z_{i_{l}}$, $1\leq l\leq 2^{\vert V \vert}$, is exactly the standard basis of $\mathbb{Z}_2^{\vert V\vert}$.
	Use $z_iS$ to denote $(\sum\limits_{v_{k} \in V_i} m_{v_k})S$, where $\varphi(z_i)=V_i$.
	By the given legal system and state, we can direct the 1-cubes at $z_ix$ by the following manner: a 1-cube is outgoing if the corresponding vertex $v$ of link$(z_ix)$ is in $z_iS$
	and it is incoming if $v\not\in z_iS$.
	
	Let $z_kx$ and $z_jx$ be the endpoints of the 1-cube $c$ corresponding to a move $v$ of $\Gamma$. And we use $z_kx\rightarrowtail z_jx$ to denote the direction of $1$-cube $c$. Suppose the symmetric difference of the vertex sets $V_k=\varphi(z_k)$ and $V_j=\varphi(z_j)$ is the vertex $v$, namely $V_k=v+V_j$. By the manner above, $z_kx\rightarrowtail z_jx$ if and only if $v\in z_kS$; $z_jx\leftarrowtail z_kx$ if and only if $v\in z_j(\Gamma-S)$. These two conditions agree with each other by the property (1) in Definition \ref{move}.
	
	Moreover, the opposite sides of each 2-cube are directed consistently. More precisely, if $z_1x,z_2x,z_3x,z_4x$ are the four vertices of a 2-cube clockwise. Without loss of generality, we assume the two adjacent vertices are $v_1$ and $v_2$ and thus, have $\varphi(z_1)=v_1+\varphi(z_2)$, $\varphi(z_3)=v_1+\varphi(z_4)$, $\varphi(z_2)=v_2+\varphi(z_3)$ and $\varphi(z_1)=v_2+\varphi(z_4)$. Therefore, $z_1S=m_{v_1}S+z_2S$, $z_3S=m_{v_1}S+z_4S$, $z_2S=m_{v_2}S+z_3S$, $z_1S=m_{v_2}S+z_4S$. By property (2) required in Definition \ref{move}, $m_{v_1}$ and $m_{v_2}$ should be $(1,0,\cdots)$ and $(0,1,\cdots)$, respectively. Then the directions of the 1-skeleton of the 2-cube can be well defined. For example, if $S$ contains vertices $v_1$ and $v_2$, and the $mz_1=(1,1,\dots)$, then we obtain the directions as shown by red arrows in Figure \ref{fig:direc}. It can be shown analogously in other cases.
	
	\begin{figure}[H]
		\scalebox{0.48}[0.48]{\includegraphics {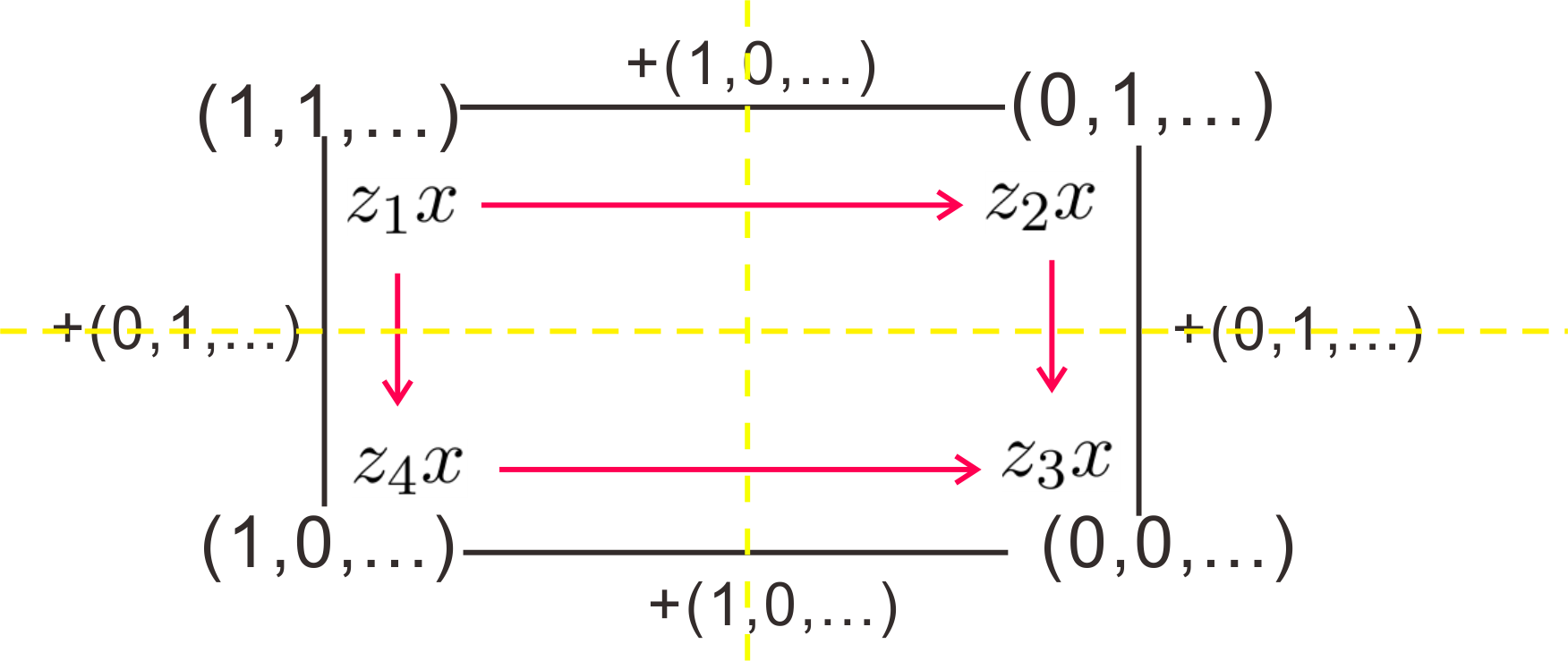}}
		\caption{Well-defined directions of the 1-skeleton of a 2-cube of the cube complex $X$.}\label{fig:direc}
	\end{figure}
\end{proof}
 By Theorem  \ref{B:finiteness} (1) and Theorem \ref{thm:win fiber}, there is a corollary as follows; the finite generation is thus  obtained:
\begin{corollary}\label{cor:win f.g.} \emph{(Jankiewicz-Norin-Wise \cite{JNW:2017})}
	Let $\Gamma$ be a simplicial graph. If $\Gamma$ admits an admissible system of moves and state, then the right-angled Coxeter group of $C(\Gamma)$ of $\Gamma$ has a finite index subgroup $G'$ such that there is an epimorphism $G'\rightarrow \mathbb{Z}$
	with finitely generated kernel is guaranteed .
\end{corollary}

\section{The fiber kernel of the real moment-angled manifold over the $24$-cell $\mathcal{P}$ with respect to the  map $\phi_{\mathcal{P}}$}
The $24$-cell $\mathcal{P}$ is one of six convex regular 4-polytopes. Its co-dimensional one faces areoctahedra; in total of $24$ There are $6$ co-dimensional one facets meeting at each vertex  and $3$ at each edge. It can be realized as a right-angled ideal hyperbolic polytope in the four-dimensional hyperbolic space $\mathbb{H}^{4}$. Let $G$ be the Coxeter group of reflections in the 3-dimensional faces of the right-angled $24$-cell. Sine the $24$-cell is self-dual, we still use $\Gamma_\mathcal{P}$ to denote the $1$-skeleton of the dual of the $24$-cell. The right-angled Coxeter group $C(\Gamma_\mathcal{P})=G$ acts properly and cocompactly on a $CAT(0)$ cube complex  $\widetilde{X}_{\mathcal{P}}$. $X_{\mathcal{P}}=G^{'}/\widetilde{X}_{\mathcal{P}}$, $G^{'}=\text{ker}~(C(\Gamma_{\mathcal{P}})\rightarrow\mathbb{Z}_2^{\vert V(\Gamma_{\mathcal{P}}) \vert})$.

We may construct a cusped real hyperbolic manifold $M(\mathcal{P})$ by the action of  $G^{'}< C(\Gamma)$ on the hyperbolic space $\mathbb{H}^4$. Such manifold is called a \emph{real moment-angle manifold} in some
 literature, for example see \cite{dj:1991}. Note that the usual real moment-angled manifold is defined over a simple polytope. For an edge-simple $n$-polytope, rather than being glued together for compact ends, the $2^n$ copies of $P$ around a vertex contribute a cusp. The smoothness are guaranteed also by \cite{Davis:2008}. In addition, $X_\mathcal{P}$ is homotopic to $M(\mathcal{P})$, thus $\pi_1(X_{\mathcal{P}})=\pi_1(M(\mathcal{P}))$. 

\subsection{The $1$-skeleton of the $24$-cell}
\begin{figure}[H]
		\scalebox{0.38}[0.38]{\includegraphics {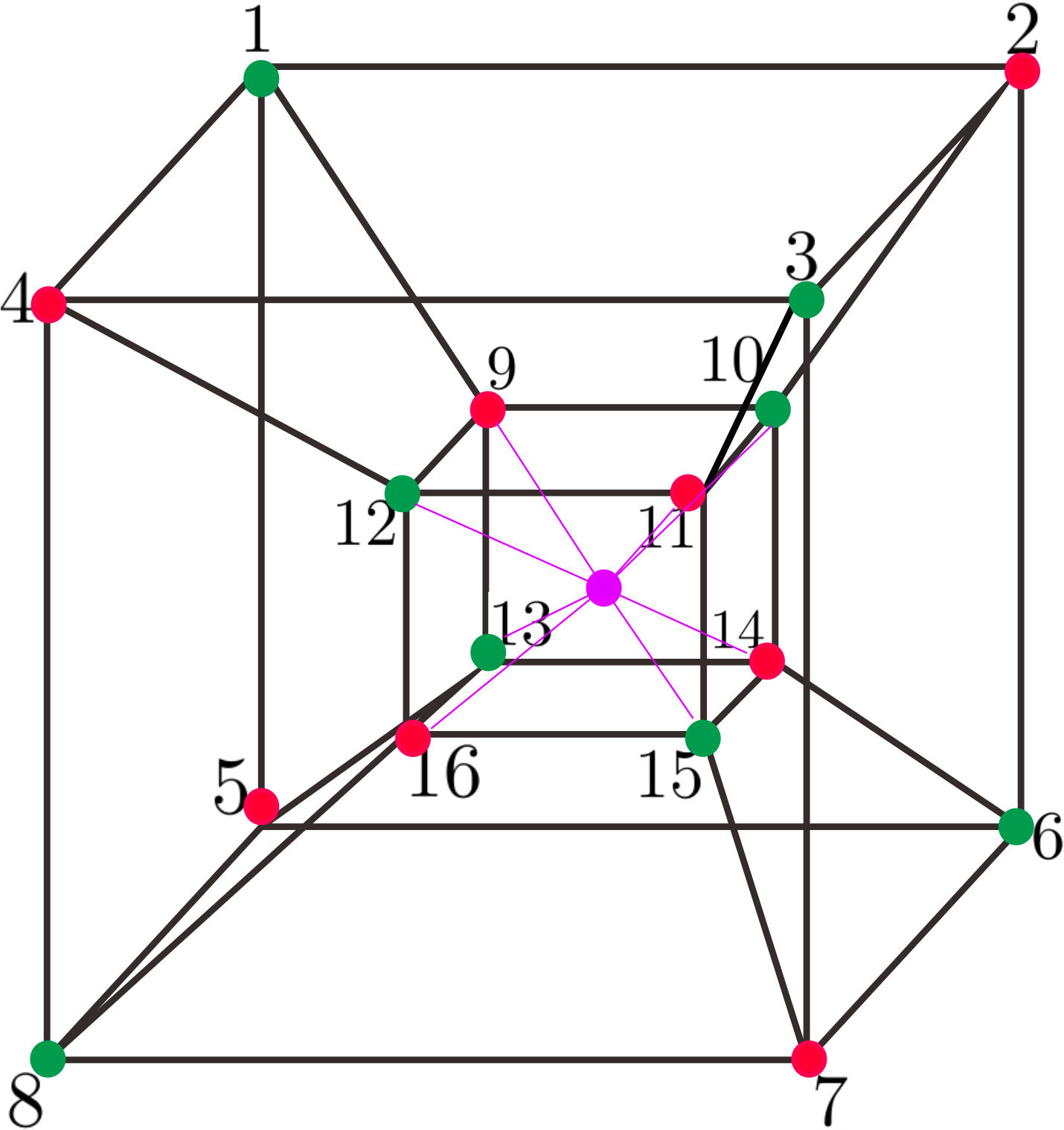}}
		\caption{The 1-skeleton of the 24-cell has 8 purple vertices corresponding to 3-dimensional faces of the 4-cube. One such vertex is illustrated here. }\label{fig:ske24}
	\end{figure}

\begin{figure}[H]
		\scalebox{0.38}[0.38]{\includegraphics {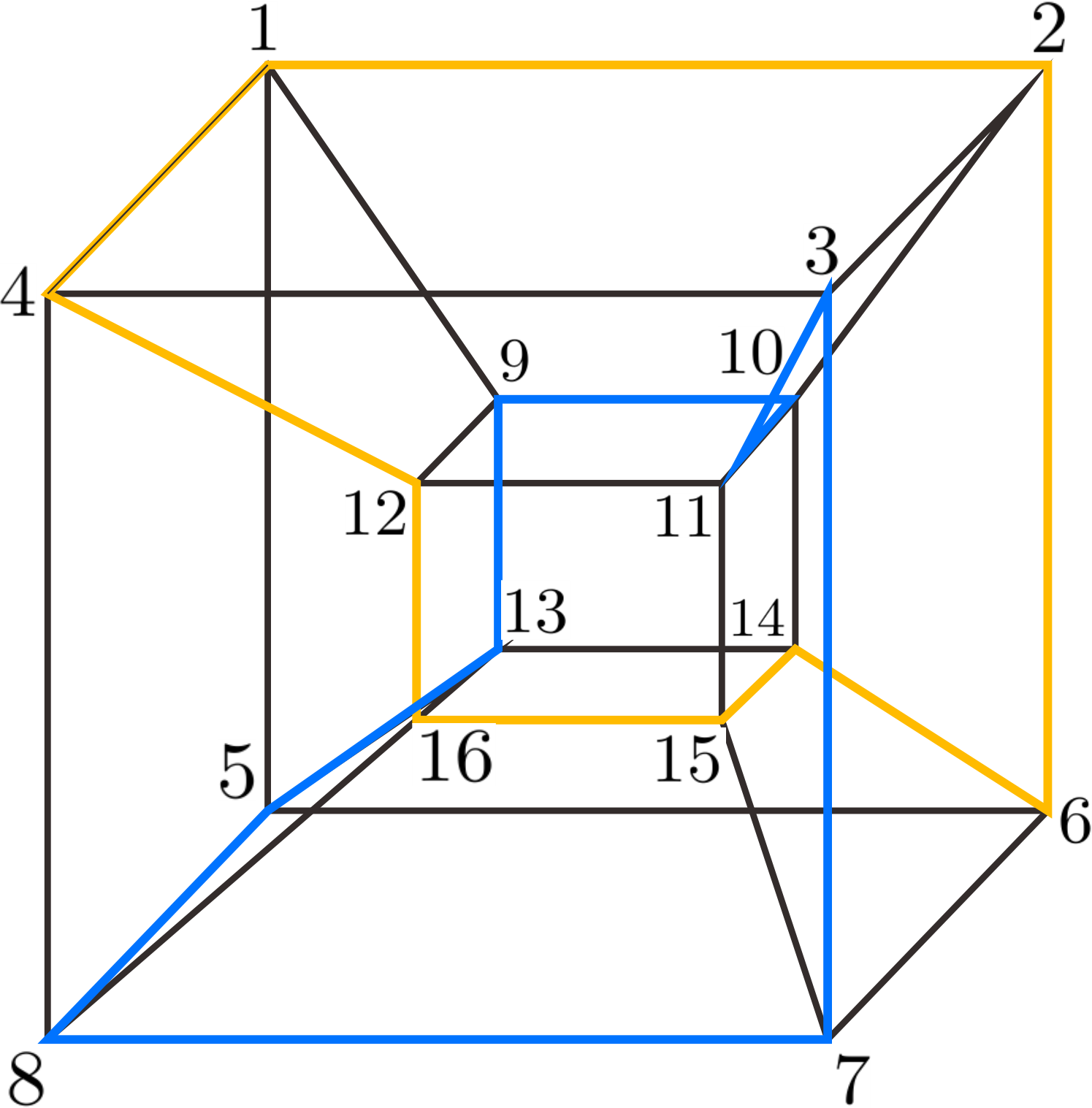}}
		\caption{A legal move system for the 1-skeleton of the
			4-cube. }\label{fig:col24}
	\end{figure}

The graph $\Gamma_{\mathcal{P}}$ can be obtained from the $1$-skeleton of the $4$-cube by adding $8$ vertices and $64$ edges. More precisely, for each of the $3$-cubes in the $4$-cube, an \emph{extra vertex} is added and connected to all vertices of its ``surrounding'' $3$-cube as shown in Figure \ref{fig:ske24}. We learn the Figure \ref{fig:ske24} and Figure \ref{fig:col24} from \cite{JNW:2017}. There are $24$ vertices in total, including $16$ vertices from the $4$-cube and $8$ extra vertices. There are totally $96$ edges, including $32$ edges from the $4$-cube, which are in black in Figure \ref{fig:ske24}, and $8$ edges for each extra vertex, which are in purple in Figure \ref{fig:ske24}. 

\subsection{Admissible system of the 1-skeleton of the $24$-cell $\Gamma_{\mathcal{P}}$.}
Jankiewicz-Norin-Wise reported in  \cite{JNW:2017} a legal move system of $\Gamma_{\mathcal{P}}$, which is obtained from a legal move system of the 1-skeleton of $4$-cube. The bipartite structure of the $4$-cube, as shown in Figure \ref{fig:col24}, naturally provides two move elements. All of the extra vertices consist of the third move. We report the admissible system  of move system $\{m_1,m_2,m_3\}$ and one starting state $S$ of  $\Gamma_{\mathcal{P}}$ constructed in \cite{JNW:2017} here:
\begin{gather}
m_1=(1,3,6,8,10,12,13,15),\notag\\
m_2=(2,4,5,7,9,11,14,16),\\
m_3=(17,18,19,20,21,22,23,24),\notag\\
S=(1,2,6,14,15,16,4,17,18,19,20).\notag
\end{gather}
And the legal orbit is as follows:
$$S=(1,2,4,6,12,14,15,16,17,18,19,20),$$
$$m_1S=(2, 3, 4, 8, 10, 13, 14, 16, 17, 18, 19, 20),$$
$$m_2S=(1, 5, 6, 7, 9, 11, 12, 15, 17, 18, 19, 20),$$
$$m_3S=(1, 2, 4, 6, 12, 14, 15, 16, 21, 22, 23, 24),$$
$$m_1m_2S=(3, 5, 7, 8, 9, 10, 11, 13, 17, 18, 19, 20),$$
$$m_1m_3S=(2, 3, 4, 8, 10, 13, 14, 16, 21, 22, 23, 24),$$
$$m_2m_3S=(1, 5, 6, 7, 9, 11, 12, 15, 21, 22, 23, 24),$$
$$m_1m_2m_3S=(3, 5, 7, 8, 9, 10, 11, 13, 21, 22, 23, 24).$$

By Corollary \ref{cor:win f.g.}, Jankiewicz-Norin-Wise \cite{JNW:2017} showed there exists an exact sequence
\begin{equation}
1\rightarrow H_{\mathcal{P}}\rightarrow \pi_1(M(\Gamma_{\mathcal{P}}))\xrightarrow{\phi_{\mathcal{P}}} \mathbb{Z}\rightarrow 1,
\end{equation}
where the fiber-kernel $H_{\mathcal{P}}$ is finitely generated.

\subsection{The fiber-kernel is not of type $FP_2$} We now show $H_{\mathcal{P}}$ in $(3.2)$
 is not finitely presented. By Theorem \ref{B:finiteness}, we have the following lemma to depict the finiteness of the fiber-kernel $H_{\mathcal{P}}$:
\begin{lemma}
Suppose $\widetilde{\phi}:\widetilde{X}\rightarrow \mathbb{R}$ is the associated Morse function with respect to the admissible system $(3.1)$ of the  given system of the 1-skeleton of the 24-cell $\Gamma_{\mathcal{P}}$. If all of the ascending links $A_i$ and descending links $D_i$ 
with respect to $\widetilde{\phi}$ satisfies

(1) $H_2(A_i; \mathbb{Z})=H_2(D_i; \mathbb{Z})=\widetilde{H_0}(A_i; \mathbb{Z})=\widetilde{H_0}(D_i; \mathbb{Z})$=0.

(2) $H_1(A_i; \mathbb{Z})\ne 0$, $ H_1(D_i; \mathbb{Z})\ne 0$.

Then the fiber kernel $H_{\mathcal{P}}$ of $\phi_{\mathcal{P}}:\pi_1(M(\Gamma_{\mathcal{P}}))\rightarrow\mathbb{Z}$ is finite-generated, but not finite-presented.
\end{lemma}

\begin{figure}[H]
		\scalebox{0.31}[0.31]{\includegraphics {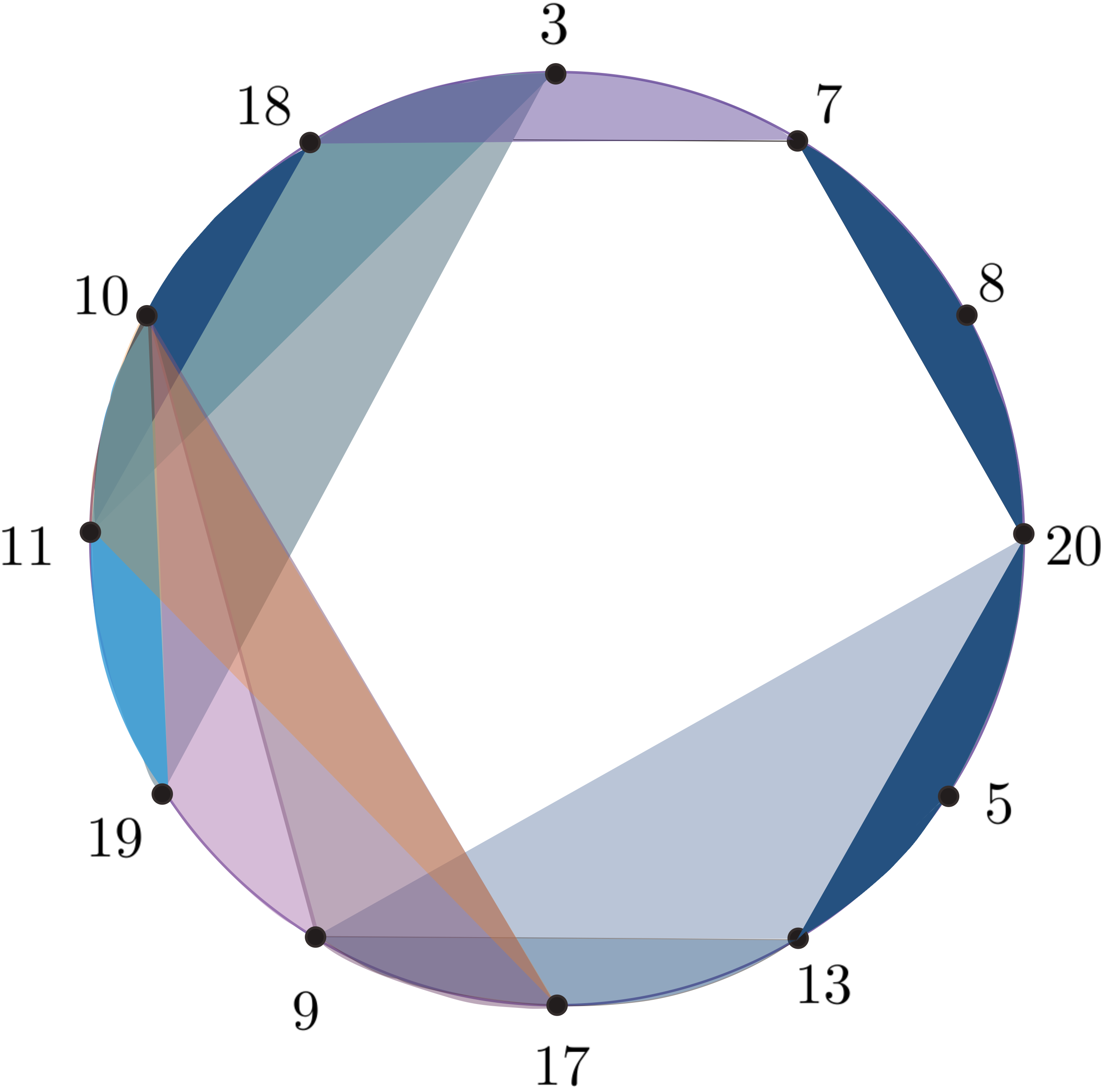}}
		\caption{Ascending link associated to $m_1m_2S$ .}\label{fig:eg}
	\end{figure}

We calculate the homologies about all of the ascending and descending links as shown in Table \ref{table:homology}. For example, the descending link spanned by the vertex set $m_1m_2S$ is as shown in Figure \ref{fig:eg}. More precisely, the 24-cell without interior is a simplicial complex;  the descending link is a full-subcomplex of it, which is  of dimension at most three and of the vetex set $m_1m_2S$. This descending link has

$\bullet$ $12$ vertices labeled by $3$, $7$, $8$, $20$, $5$, $13$, $17$, $9$, $19$, $11$, $10$, $18$;

$\bullet$ $24$  edges labeled by $(3,7)$, $(3,19)$, $(3,11)$, $(7,8)$, $(7,20)$, $(7,18)$, $(8,20)$, $(20,5)$, $(20,13)$, $(20,9)$,   $(5,13)$, $(13,17)$, $(13,9)$, $(17,9)$, $(17,11)$, $(17,10)$,  $(9,19)$, $(9,10)$, $(19,11)$, $(19,10)$, $(11,10)$, $(11,18)$, $(10,18)$, $(18,3)$;

$\bullet$ $12$  triangles labeled by $(3,7,18)$, $(3,19,11)$, $(3,11,18)$,  $(7,8,20)$, $(20,5,13)$, $(20,13,9)$, $(13,17,9)$, $(17,9,10)$, $(17,11,10)$, $(9,19,10)$, $(19,11,10)$, $(11,10,18)$.

$\bullet$ no tetrahedron.

This descending link  is homotopic to $\mathbb{S}^1$.  Homologies of the ascending and descending links with respect to $\phi_{\mathcal{P}}$ 
are as shown in Table \ref{table:homology}.  In particular, the descending link spanned by the vertex set $m_3S$ is homotopic to $\mathbb{S}^1\vee\mathbb{S}^1$.
\begin{table}[H]
	\small
	\begin{tabular}{|c|c|c|c|c|c|c|c|c|}
		\hline
		State $S_i$ that spans the $D_i$&$S$&$m_1S$&$m_2S$&$m_3S$&$m_1m_2S$&$m_1m_3S$&$m_2m_3S$&$m_1m_2m_3S$\\
		\hline	
		 $\widetilde{H_0}(D_i; \mathbb{Z})$&0&0&0&0&0&0&0&0\\
		 \hline
		$H_1(D_i; \mathbb{Z})$ &$\mathbb{Z}$&$\mathbb{Z}$&$\mathbb{Z}$&$\mathbb{Z}\oplus\mathbb{Z}$ &$\mathbb{Z}$ &$\mathbb{Z}$&$\mathbb{Z}$&$\mathbb{Z}$\\
		\hline
		$H_2(D_i; \mathbb{Z})$&0&0&0&0&0&0&0&0\\
		\hline
		\hline
		State $S_i$ that spans the  $A_i$&$m_1m_2m_3S$&$m_2m_3S$&$m_1m_3S$&$m_1m_2S$&$m_3S$&$m_2S$&$m_1S$&$S$\\
		\hline	
		$\widetilde{H_0}(A_i; \mathbb{Z})$&0&0&0&0&0&0&0&0\\
		\hline
		$H_1(A_i; \mathbb{Z})$ &$\mathbb{Z}$&$\mathbb{Z}$&$\mathbb{Z}$&$\mathbb{Z}$ &$\mathbb{Z}\oplus\mathbb{Z}$ &$\mathbb{Z}$&$\mathbb{Z}$&$\mathbb{Z}$\\
		\hline
		$H_2(A_i; \mathbb{Z})$&0&0&0&0&0&0&0&0\\
		\hline
	\end{tabular}
	\caption{Homologies of the ascending and descending links with respect to $\phi_{\mathcal{P}}$ at each $0$-cube of $X_\Gamma$}
	\label{table:homology}
\end{table}

By Table \ref{table:homology}, we complete the proof of Theorem \ref{main} for the real moment angle manifold over the 24-cell $\mathcal{P}$.

\section{The fiber kernel of the real moment-angled manifold over the $120$-cell $\mathcal{E}$ with respect to the map $\phi_{\mathcal{E}}$}
The 120-cell $\mathcal{E}$ is a notable convex regular 4-polytope that admits an embedding
in the four-dimensional hyperbolic space $\mathbb{H}^{4}$ as a right-angled polytope. The dual of the 120-cell is the 600-cell that has 600 co-dimensional one faces of tetrahedra and the link of each vertex is an icosahedron. We use $\Gamma_{\mathcal{E}}$ to denote the graph of 1-skeleton of 600-cell.

\subsection{The $1$-skeleton of the $600$-cell}

\begin{figure}[H]
		\scalebox{0.5}[0.5]{\includegraphics {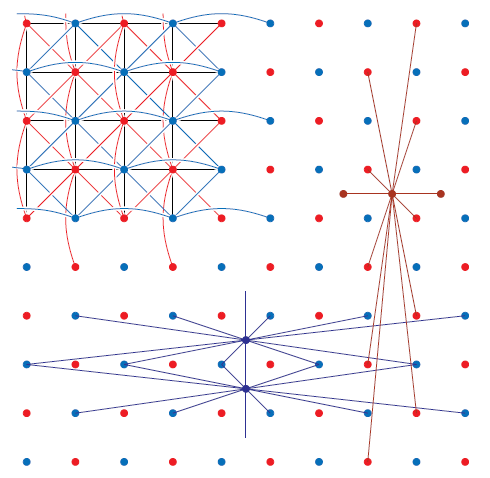}}
		\caption{Illustration of the 1-skeleton of the 600-cell. Red and blue
			vertices are even and odd respectively. One even hovering vertex and two
			consecutive odd hovering vertices are illustrated. }\label{fig:600ske}
	\end{figure}

The combinatoric knowledge about the 600-cell is collected from Jankiewicz-Norin-Wise's paper \cite{JNW:2017} and from the web \cite{600-cell}. We report them here for the convenience of readers. We learn the Figure \ref{fig:600ske} and Figure \ref{fig:600mov} from \cite{JNW:2017}.

Let's prepare a torus represented by a $10\times 10$ grid ahead. 
The combinatorics of the 1-skeleton of the 600-cell can be depicted as follows:

$\bullet$ There are in total of 120 vertices:

\begin{enumerate}
	\item 100 of them are given by the $10\times 10$ grid and named \emph{ordinary vertices}. These vertices are bi-partied into \emph{even vertices} and \emph{odd vertices}, which are in red and blue, respectively, in Figure \ref{fig:600ske};	
	\item The other $20$ vertices named \emph{hovering vertices}. Ten of them are in consecutive pairs of vertical cycles, which are named \emph{even hovering vertices} and colored by red in Figure \ref{fig:600ske}. Whereas, the other ten are in consecutive pairs of horizontal cycles, which are named \emph{odd hovering vertices} and colored by blue in Figure \ref{fig:600ske}. 
\end{enumerate} 

$\bullet$ There are in total of 720 edges:
\begin{enumerate}
	\item 200 of them are given by the grid edges;
	\item 100 of them are given by joining the ordinary vertices even-to-even and odd-to-odd as shown in Figure \ref{fig:600ske} ;
	\item 100 of them are given by adding two diagonals in each  square of the grid as shown in Figure \ref{fig:600ske} 
	\item 100 of them are given by attaching each even hovering vertex to all ten even vertices of the two vertical cycles that the even hovering vertex is in between as shown is Figure \ref{fig:600ske}. Whereas 100 of them are given by attaching each odd hovering vertex to all ten odd vertices of the two horizontal cycles that the odd hovering vertex is in between as shown is Figure \ref{fig:600ske}.
	\item 20 of them are given by adding edges between consecutive even hovering vertices and consecutive odd hovering vertices
\end{enumerate}

\subsection{Admissible system for the 1-skeleton of the $600$-cell $\Gamma_{\mathcal{E}}$.} \label{system120}
By \cite{JNW:2017}, there is a system of moves for the 600-cell's 1-skeleton $\Gamma_{\mathcal{E}}$. The vertices of the grid-vertices of $\Gamma_{\mathcal{E}}$ are labelled by $\{0,1,2,3,4,5,6,7,8,9\}$ as on the right in Figure \ref{fig:600mov}:

	\begin{figure}[H]
		\scalebox{0.46}[0.46]{\includegraphics {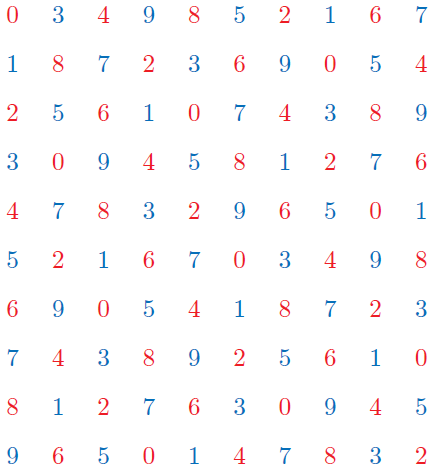}}
		\caption{a legal system over the 1-skeleton of the 600-cell. Moves correspond to all vertices labelled with the same number. There are also moves corresponding to the hovering vertices which are not illustrated in the figure.}\label{fig:600mov}
	\end{figure}

The even hovering vertices are labelled by $0', 2', 4', 6', 8', 0', 2', 4', 6', 8'$ consecutively. The odd hovering vertices are labelled by $1', 3', 5', 7', 9', 1', 3', 5', 7', 9'$ consecutively.

The start state $S_0$ consists of all $60$ vertices in alternate horizontal circles of the torus, together with a choice of alternate vertices in both even and odd hovering cycle as shown in Figure \ref{fig:state600}.

Consider the move system consisting of 20 moves $m_{1}, m_{2}, \cdots, m_{20}$, which  corresponding to the distinct labels $\{0,1,2,3,4,5,6,7,8,9, 0',1',2',3',4',5',6',7',8',9'\}$ for the vertex set of $\Gamma_{\mathcal{E}}$. By Corollary \ref{cor:win f.g.}, Jankiewicz-Norin-Wise  \cite{JNW:2017} reported an exact sequence
\begin{equation}
1\rightarrow H_{\mathcal{E}} \rightarrow \pi_1(M(\mathcal{E}))\xrightarrow{\phi_{\mathcal{E}}} \mathbb{Z}\rightarrow 1,
\end{equation}

\noindent where $H_{\mathcal{E}}$ is finitely generated;  $\phi_{\mathcal{E}}$ is induced from a  Bestvina-Brady's  Morse function and is constructed by Jankiewicz-Norin-Wise in \cite{JNW:2017}.  Here $M_{\mathcal{E}}$ is the real moment-angle manifold over the 120-cell $\mathcal{E}$ and $\pi_1(M(\mathcal{E}))$, by definition, algebraically fibers. 

The right-angled Coxeter group $C(\Gamma_\mathcal{E})=G$ acts properly and cocompactly on a $CAT(0)$ cube complex  $\widetilde{X}_{\mathcal{E}}$. $X_{\mathcal{E}}=G^{'}/\widetilde{X}_{\mathcal{E}}$, $G^{'}=\text{ker}~(C(\Gamma_{\mathcal{E}})\rightarrow\mathbb{Z}_2^{\vert V(\Gamma_{\mathcal{E}}) \vert})$. Therefore, $X_{\mathcal{E}}$ is the 4-dimensional cube complex homotopic to $M(\mathcal{E})$. The set  $\{m_1,m_2,\cdots, m_{20}\}$ is a move system $\mathcal{M}$ of the 1-skeleton $\Gamma_{\mathcal{E}}$ of the 120-cell $\mathcal{E}$ . The sets of moves generate a subgroup of $\mathbb{Z}_2^{\vert V(\Gamma_{\mathcal{E}})\vert}$, where $V(\Gamma_{\mathcal{E}})$ is the vertex set of $\Gamma_{\mathcal{E}}$. Denote the ascending and descending links by $A_i$ and $D_i$, respectively. Although $\vert \langle m_1,m_2,\cdots, m_{20}\rangle\vert=2^{20}$, we do not need to calculate all of the  pairs of reduced homology groups of the ascending link $A_i$ and descending link $D_i$ with respect to the vertex $v_i$ of the cube complex $\widetilde{X}_{\mathcal{E}}$. By referring details of Theorem 4.7 in Brady's paper \cite{Brady:1999}, it turns out that we only need to ensure the following two conditions to claim that the fiber-kernel  $H_{\mathcal{E}}$ in $(4.1)$ is not $FP_2$:

(a) $\widetilde{H_i}(A_j; \mathbb{Z})=0,\widetilde{H_i}(D_j; \mathbb{Z})=0$, $i=0,2$, $j=1,2,...,2^{20}.$
	
(b) There exists  some descending link, for example $D_1$, such that $H_1(D_1; \mathbb{Z})\ne 0$.

 \subsection{The fiber-kernel is not of type $FP_2$}
Now we are going to verify  conditions (a) and (b)  proposed above.

\noindent \textbf{Proof of condition (a):}

We use the notation $\phi_{\mathcal{E}}$, $X$, $A_i$, $D_j$ as claimed above. By \cite{JNW:2017}, the ascending links $A_i$ and descending links $D_j$ are always connected. We need to show that $H_2(A_i; \mathbb{Z})=H_2(D_j; \mathbb{Z})=0$.

\begin{figure}[H]
		\scalebox{0.85}[0.85]{\includegraphics {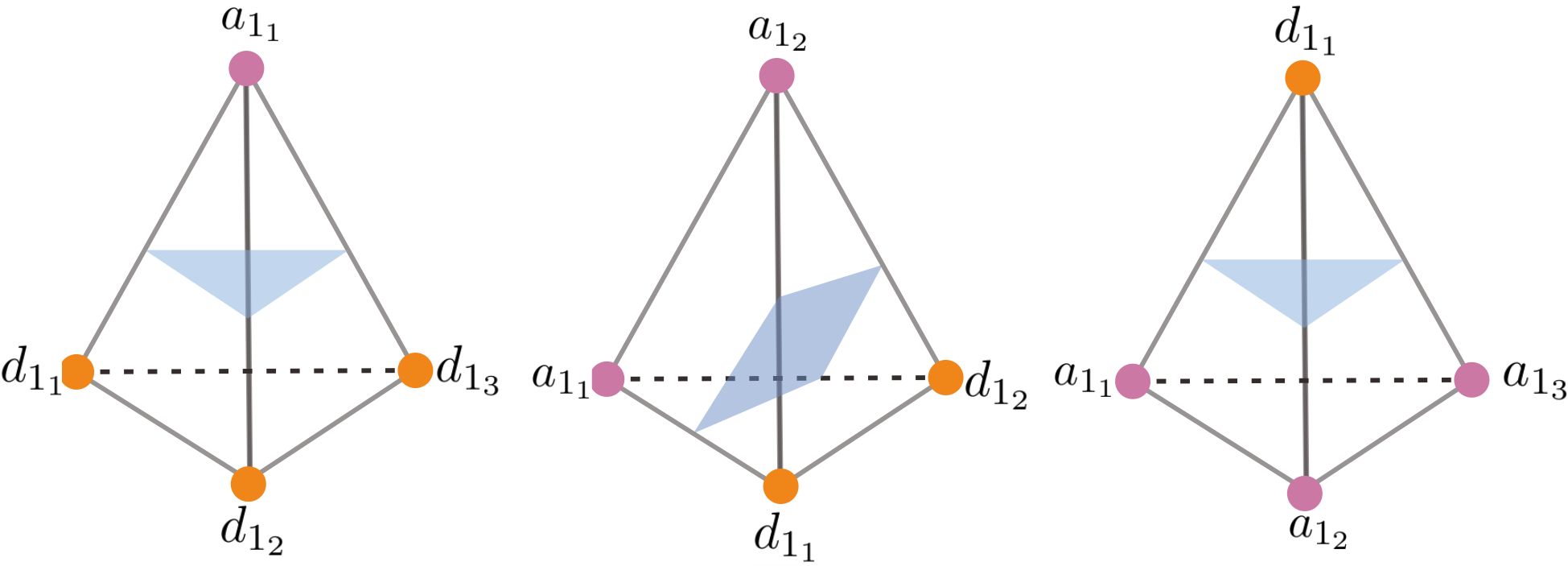}}
		\caption{There is always a canonical hyperplane in a tetrahedron that separates the vertices  belonging to $A_i$, denoted by $a_{i_{j}}$, and $D_i$, denoted by $d_{i_{k}}$. }\label{fig:surface}
	\end{figure}

The vertex linking of a vertex of the cube complex $\widetilde{X}_{\mathcal{E}}$ is always a 600-cell, which is homeomorphic to $\mathbb{S}^3$. We still use the notation $A_i$ and $D_i$ to denote the parts of $\mathbb{S}^3$ that $A_i$ and $D_i$ are projected onto. We claim that $A_i$ and $D_i$ can be separated by a connected closed orientable surface $F_i$ in $\mathbb{S}^3$.
First, the boundary of 600-cell consists of 600 tetrahedra. In one tetrahedron, there is a canonical way to construct a hyperplane to separate the vertices in $A_i$ and $D_i$  as shown in Figure \ref{fig:surface}. All of these hyperplanes can be glued together to form a closed surface $F_i$ in $\mathbb{S}^3$ that separating  $A_i$ and $D_i$. Moreover, $F_i$ is guaranteed to be orientable because every closed surface in $\mathbb{S}^3$ is orientable. Denote the parts of $\mathbb{S}^3$ that contains $A_i$ and $D_i$ by $\widetilde{A}_i$ and $\widetilde{D}_i$, respectively. Then $\widetilde{A}_i\cap \widetilde{D}_i=F_i$, $\widetilde{A}_i\cup \widetilde{D}_i=\mathbb{S}^3$, $A_i$ and $D_i$ are strong deformations of $\widetilde{A}_i$ and $\widetilde{D}_i$, respectively.

\begin{figure}[H]
		\scalebox{0.4}[0.4]{\includegraphics {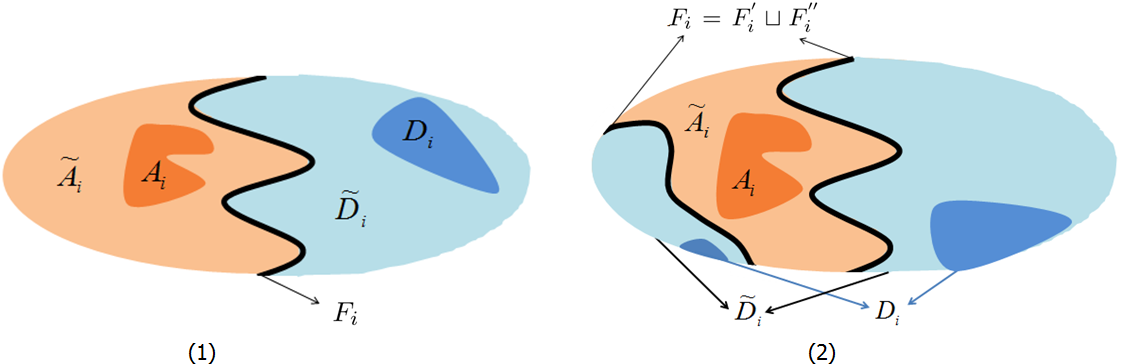}}
		\caption{}\label{fig:f}
	\end{figure}

  Next, if the separating surface is not connected, in other words, $F_i=F_i^{'}\sqcup F_i^{''}$, since every closed surface in $\mathbb{S}^3$ is separating, then $\widetilde{A}_i$ and $\widetilde{D}_i$ are not connected either as shown in Figure \ref{fig:f}.

 That contradicts with the assumption that  $A_i$ and $D_i$ are connected. Thus $\widetilde{A}_i\cap \widetilde{D}_i=F_i$ is a closed connected orientable surface.

	Now, for each $i$, we have a Mayer-Vietoris sequence
	$$0\rightarrow H_3(\widetilde{A}_i\cup \widetilde{D}_i; \mathbb{Z})\xrightarrow{f_*}H_2(\widetilde{A}_i\cap\widetilde{D_i}; \mathbb{Z})\rightarrow H_2(\widetilde{A}_i; \mathbb{Z})\oplus H_2(\widetilde{D}_i; \mathbb{Z})\rightarrow 0,$$	
\noindent where $f_*=\Delta^3\circ \tau \circ i_*$ is an isomorphism induced by the long exact sequence
	$$\cdots\rightarrow H_3(\widetilde{A}_i\cup \widetilde{D}_i; \mathbb{Z})\xrightarrow{i_*}  H_3(\widetilde{A}_i\cup \widetilde{D}_i,\widetilde{D}_i; \mathbb{Z})\xrightarrow{\tau} H_3(\widetilde{A}_i,\partial\widetilde{A}_i=\partial\widetilde{D}_i; \mathbb{Z})\xrightarrow{\Delta^3}    H_2(\widetilde{A}_i\cap\widetilde{D}_{i}; \mathbb{Z})\rightarrow \cdots,$$
	\noindent where $\tau$ is induced by exclusion.  Therefore, $H_2(A_i; \mathbb{Z})=H_2(\widetilde{A}_i; \mathbb{Z})=0$, $H_2(D_i; \mathbb{Z})=H_2(\widetilde{D}_i; \mathbb{Z})=0$.
	
\noindent \textbf{Proof of condition (b):}
    In the following, we calculate that $H_1(D_1; \mathbb{Z})\ne 0$, where descending link $D_1$ is the full-subcomplex of the 600-cell that defined on the starting state $S_0$.

	\begin{figure}[H]
		\scalebox{0.5}[0.7]{\includegraphics {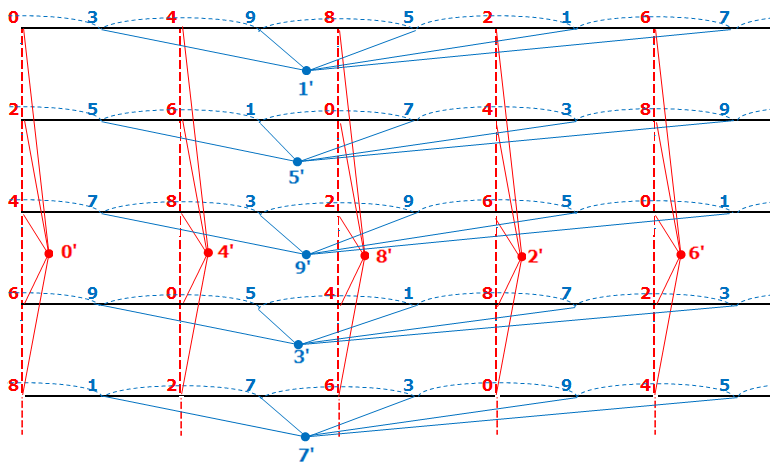}}
		\caption{The starting state $S_0$ of the 1-skeleton of 600-cell}\label{fig:state600}
	\end{figure}

Because the descending link is connected and we already have $H_2(D_1; \mathbb{Z})=0$, the Euler characteristic $\chi(D_1)=\beta_0(D_1)-\beta_1(D_1)+\beta_2(D_1)-\beta_3(D_1)=1-\beta_1(D_1)$.

On the other hand, we can calculate $\chi(D_1)$ by counting simplices of dimension $0$, $1$, $2$ and $3$ in the simplicial complex $D_1$. As shown in Figure \ref{fig:state600}, there are in total of 60 vertices in $D_1$, then $e_0=|S_{0}|=60$. Moreover, 150 edges can be found in $D_1$, including 50 edges of type 	\scalebox{0.4}[0.5]{\includegraphics {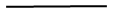}}, 25 edges of type \scalebox{0.4}[0.5]{\includegraphics {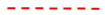}}, 25 edges of type 	 \scalebox{0.4}[0.5]{\includegraphics {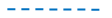}}, 25 edges of type 	 \scalebox{0.4}[0.5]{\includegraphics {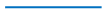}}, 25 edges of type 	\scalebox{0.4}[0.5]{\includegraphics {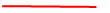}}. Namely $e_1=150$. There are 75 triangles in total in $D_1$, including 25 triangles of type 	\scalebox{0.4}[0.5]{\includegraphics {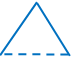}}, 25 triangles of type 	\scalebox{0.4}[0.5]{\includegraphics {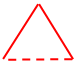}}, 25 triangles of type 	\scalebox{0.4}[0.5]{\includegraphics {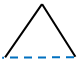}}, so  $e_2=75$. No tetrehedra can be obtained in $D_1$ and thus $e_3=0$. So we have $\chi(D_1)=60-150+75-0=15=\beta_0(D_1)-\beta_1(D_1)$. Therefore $\beta_1(D_1)=16 \ne 0$.

Now, We complete the proof of Theorem \ref{main} for the real moment angle manifold over the 120-cell $\mathcal{E}$.


\begin{thebibliography}{10}
	
\bibitem{Agol:2008}	
I. Agol.
\newblock  Criteria for virtual fibering.
\newblock   J. Topol., 1 (2008), no. 2, 269--284.
	
\bibitem{Agol:2013}	
I. Agol.
\newblock  The virtual Haken conjecture. With an appendix by I. Agol, D. Groves and J. Manning.
\newblock   Doc. Math., 18 (2013), 1045--1087.

\bibitem{as:2019}	
I. Agol and M. Stover.
\newblock  Congruence RFRS towers.
\newblock   ArXiv:1912.10283.

\bibitem{Bieri:1976}
R. Bieri.
\newblock Homological dimension of discrete groups.
\newblock Queen Mary College Mathematics Notes, 1976.


\bibitem{BB:97}
M. Bestvina and N. Brady.
\newblock Morse theory and finiteness properties of groups.
\newblock Invent. Math., 129 (1997), no.3, 445--470.


\bibitem{Bowditch_Mess:1994}	
B. H. Bowditch and G. Mess.
\newblock A 4-dimensional Kleinian group.
\newblock Trans. Amer. Math. Soc., 344 (1994), no. 1, 391--405.

\bibitem{Brady:1999}
N. Brady.
\newblock Branched coverings of cubical complexes and subgroups of hyperbolic groups.
\newblock J. Lond. Math. Soc., 60 (1999), no. 2,  461--480.

\bibitem{Brown:1982}
K. Brown.
\newblock Cohomology of groups.
\newblock GTM 87, Springer-verlag, New York, 1982.



\bibitem{Davis:2008}	
M. Davis.
\newblock The geometry and topology of Coxeter groups.
\newblock Volume 32 of London Mathmetical Society Monographs Series, Princeton University Press, Princeton, NJ, 2008.


\bibitem{dj:1991}
M. Davis and T. Januszkiewicz.
\newblock  Convex polytopes, Coxeter orbifolds and torus actions.
\newblock   Duke Math. J. 62 (1991), 417--451.

\bibitem{FV}
S. Friedl and S. Viduss.
\newblock On virtual properties of K\"{a}hler groups.
\newblock ArXiv:1704.07041, to appear in Nagoya Math. J.


\bibitem{fv:2019}	
S. Friedl and S. Vidussi.
\newblock  BNS invariants and algebraic fibrations of group extensions.
\newblock   ArXiv:1912.10524.


\bibitem{JNW:2017}
K. Jankiewicz, S. Norin and D. Wise.
\newblock Virtually fibering right-angled coxeter groups.
\newblock Journal of the Institute of Mathematics of Jussieu, published online by Cambridge University Press: 23 August 2019.DOI: https://doi.org/10.1017/S1474748019000422





\bibitem{Kapo:1998}
M. Kapovich.
\newblock On normal subgroups of the fundamental groups of complex surfaces.
\newblock Preprint.

\bibitem{Kapo:2013}
M. Kapovich.
\newblock Noncoherence of arithmetic hyperbolic lattices.
\newblock Geom. Topol., 17 (2013), no. 1, 39--71.


\bibitem{Kielak:2018:RFRS}
D. Kielak.
\newblock Residually finite rationally-solvable groups and virtual fibering.
\newblock Arxiv:1809.09386v1, to appear in Jour. Amer. Math. Soc.


\bibitem{Krop:2018}
R. Kropholler.
\newblock Hyperbolic groups with almost finitely presented.
\newblock ArXiv:1809.10594.

\bibitem{Liu:2013}	
Y. Liu.
\newblock  Virtual cubulation of nonpositively curved graph manifolds.
\newblock   J. Topol., 6 (2013), no. 4, 793--822.


\bibitem{Stallings:1962}
J. Stallings.
\newblock On fibering certain 3-manifolds.
\newblock 1962 Topology of 3-manifolds and related topics (Proc. The Univ. of Georgia Institute, 1961) pp. 95--100 Prentice-Hall, Englewood Cliffs, N.J.

\bibitem{Stover:2018}	
M. Stover.
\newblock Cusp and $b_1$ growth for ball quotients and maps onto $\mathbb{Z}$ with finitely generated kernel.
\newblock To appear in Indiana Univ. Math. J.


\bibitem{Wise:2018}	
D. T. Wise.
\newblock  The structure of groups with a quasiconvex hierarchy.
\newblock  
To appear in Ann. Math. Studies.


\bibitem{600-cell}	https://en.wikipedia.org/wiki/600-cell.




	
\end{thebibliography}
\end{document}